\documentclass{amsart}
\usepackage{amssymb,latexsym,graphicx,psfrag,epstopdf,hyperref,color,pinlabel}

\usepackage[all,arc,dvips]{xy}

\newcommand{\executeiffilenewer}[3]{%
\ifnum\pdfstrcmp{\pdffilemoddate{#1}}%
{\pdffilemoddate{#2}}>0%
{\immediate\write18{#3}}\fi%
}

\newcommand{\p}{\pi}
\def\frac#1#2{{#1\over#2}}

\newcommand{\QQ}{\mathbb{Q}}
 
\newcommand{\RR}{\mathbb{R}}

\newcommand{\VV}{\mathcal{V}}

\def\vol{\mathrm{Vol}}

\def\area{\mathrm{Area}}

\def\th{\theta}

\def\Vol{\mathrm{Vol}}
\def\area{\mathrm{Area}}

\def\HH{\mathbb{H}}

\def\MC{\mathcal{C}}

\def\MO{\mathcal{O}}
\def\PP{\mathcal{P}}

\def\QQ{\mathcal{Q}}

\def\co{\colon\thinspace}

\def\MC{\mathcal{C}}

\def\MO{\mathcal{O}}
\def\PP{\mathcal{P}}

\def\QQ{\mathcal{Q}}

\def\MC{\mathcal{C}}

\def\MO{\mathcal{O}}

\newtheorem{proposition}{Proposition}[section]
\newtheorem{theorem}[proposition]{Theorem}

\newtheorem{lemma}[proposition]{Lemma}

\newtheorem*{thma}{Theorem \ref{T:SmallestLarge}}

\theoremstyle{remark}

\theoremstyle{remark}

\newtheorem*{acknowledgments}{Acknowledgments}

\numberwithin{equation}{section}

\def\Area{\mbox{\rm{Area}}}
\def\Vol{\mbox{\rm{Vol}}}

\def\co{\colon\thinspace}

\begin{document}
\Large
\title{The smallest {H}aken hyperbolic polyhedra}
\author[]{Christopher~K.~Atkinson}
\address{Department of Mathematics, Temple University, Philadelphia, PA
19106, USA}
\email{ckatkin@temple.edu}
\author[]{Shawn~Rafalski}
\address{Department of Mathematics and Computer Science, Fairfield University, Fairfield, CT 06824, USA}
\email{srafalski@fairfield.edu}
\date{\noindent August 2011. 
}

\begin{abstract}
We determine the lowest volume hyperbolic Coxeter polyhedron whose 
corresponding hyperbolic polyhedral $3$--orbifold contains an essential 
$2$--suborbifold, up to a canonical decomposition along essential 
hyperbolic triangle $2$--suborbifolds.   
\end{abstract}
\maketitle

\section{Introduction}\label{S:Intro}
The organization of the volumes of hyperbolic $3$--manifolds and $3$--orbifolds 
is ongoing. Gabai, Meyerhoff and Milley have identified the Weeks--Fomenko--Matveev 
manifold as the lowest volume hyperbolic $3$--manifold 
\cite{GabMeyMill1, GabMeyMill2}, and 
Gehring, Marshall and Martin have identified the lowest volume 
hyperbolic $3$--orbifolds \cite{GehringMartin, MarshallMartin}.
Restricting to the case of orbifolds, a natural class to consider, 
from the standpoint of volume organization, is that of polyhedral 
$3$--orbifolds, i.e., the orientable $3$--orbifolds that correspond to tilings of 
hyperbolic $3$--space $\mathbb{H}^{3}$ by finite volume Coxeter polyhedra. 
For any hyperbolic $3$--orbifold that is diffeomorphic to the interior of a compact 
orbifold, there is a canonical decomposition, due to Dunbar, 
along totally geodesic hyperbolic turnovers (which are  
quotients of $\mathbb{H}^{2}$ by hyperbolic triangle groups) into components 
that either contain 
an embedded, essential $2$--suborbifold or contain no 
embedded, essential $2$--suborbifolds.
Moreover, because the decomposing turnovers are totally geodesic, 
this decomposition is volume additive, and so it is natural to consider the 
lowest volume hyperbolic polyhedral $3$--orbifolds (or their associated 
polyhedra) that either contain or do not contain essential $2$--suborbifolds, 
up to the Dunbar decomposition. 
The lowest volume polyhedral orbifolds in the latter case (called \emph{small} orbifolds) 
have been identified by Rafalski \cite{rasmallski}. 
This paper addresses the former case, 
or what are called \emph{Haken} polyhedra.

We prove the following theorem:

\begin{theorem}\label{T:SmallestLarge}
The smallest volume Haken hyperbolic Coxeter polyhedron is the Lambert
cube $\MC$.
\end{theorem}

\begin{figure}
		\centering
		\def\svgwidth{2.5in}
	    	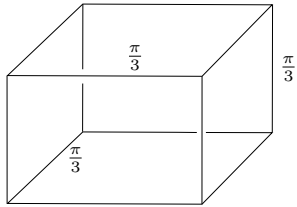
	      \caption{The smallest hyperbolic Coxeter polyhedron 
		whose associated polyhedral $3$--orbifold 
		contains an essential $2$--orbifold but contains no 
		essential triangular $2$--orbifolds. 
		All unlabeled edges have dihedral angle $\pi/2$.}		
		\label{F:LambertCube}
\end{figure}

The Lambert cube is combinatorially a cube with all but three of its 
dihedral angles equal to $\pi/2$. The remaining three dihedral angles 
equal $\pi/3$, as in Figure~\ref{F:LambertCube}. 
The Lambert cube is an example of a \emph{hyperbolic 
Coxeter $n$--prism}. For each $n\geq5$, there is a (Haken) $n$--prism that is
conjectured to be the lowest volume hyperbolic polyhedron with $2n$ vertices. 

The main idea of the proof of Theorem~\ref{T:SmallestLarge} is to use various
volume bounds for hyperbolic polyhedra to restrict the combinatorial type of
a polyhedron with volume less than the volume of the Lambert cube.  The
organization of the paper is as follows:  In Section~\ref{s:background} we
discuss the background on hyperbolic polyhedra and orbifolds in order to
understand the main result.  In Section~\ref{s:volumes} we prove a volume
bound for certain types of hyperbolic polyhedra.  This volume bound will be
used to prove the main theorem.  In Section~\ref{s:proof}, we prove
Theorem~\ref{T:SmallestLarge}.

\begin{acknowledgments}
The authors wish to thank Dave Futer for helpful feedback.
\end{acknowledgments}

\section{Background and Definitions}\label{s:background}

In this section we recall the basics of hyperbolic polyhedra and relevant
definitions for orbifolds.

An \emph{abstract polyhedron} is a cellulation of $S^2$ realizable as a
convex Euclidean polyhedron.  A theorem of Steinitz \cite{steinitz} says that realizability
as a convex Euclidean polyhedron is equivalent to the $1$--skeleton of the
cellulation being a $3$--connected planar graph 
(a graph is \emph{$3$--connected} if the removal of any two vertices 
along with their incident open edges leaves the complement connected).
A \emph{labeled abstract polyhedron} $(P,\Theta)$ is an abstract polyhedron $P$
along with a function $\Theta \co \text{Edges}(P) \to (0,\pi/2]$ labeling the
edges by real numbers which should be thought of as dihedral angles.  A
labeled abstract polyhedron $(P, \Theta)$ is \emph{realizable} as a
hyperbolic polyhedron $\PP$ if there is a label preserving graph isomorphism
between $(P, \Theta)$ and the $1$--skeleton of $\PP$ labeled by dihedral
angles.

A \emph{$k$--circuit} is a simple closed curve of $k$ edges in $P^*$, the
planar dual to $P$.  A \emph{prismatic $k$--circuit} is a $k$--circuit such
that no two edges lie in a common face of $P^*$.

Andreev's theorem characterizes non--obtuse hyperbolic polyhedra in terms of
combinatorial conditions on their $1$--skeleta.  The proof of
Theorem~\ref{T:SmallestLarge} uses Andreev's theorem to restrict the
combinatorics of certain polyhedra \cite{andreev1, andreev2}.

\begin{theorem}[Andreev's theorem]\label{and}
A non--obtuse labeled abstract polyhedron $(P, \Theta)$ that has more than
$4$ vertices is realizable as a finite
volume hyperbolic polyhedron if and only if the following hold:
\begin{enumerate}
\item Each vertex meets $3$ or $4$ edges.
\item If $e_i,$ $e_j,$ and $e_k$ share a vertex then
$\Theta(e_i)+\Theta(e_j)+\Theta(e_k) \geq \pi$.
\item If $e_i,$ $e_j,$ $e_k,$ and $e_l$ share a vertex then
$\Theta(e_i)+\Theta(e_j)+\Theta(e_k)+\Theta(e_l) = 2 \pi$.
\item If $e_i,$ $e_j,$ and $e_k$ form a prismatic $3$--circuit, then 
$\Theta(e_i)+\Theta(e_j)+\Theta(e_k) < \pi.$
\item If $e_i,$ $e_j,$ $e_k,$ and $e_l$ form a prismatic $4$--circuit, then 
$\Theta(e_i)+\Theta(e_j)+\Theta(e_k)+ \Theta(e_l) < 2\pi.$
\item If $P$ has the combinatorial type of a triangular prism with edges
$e_i,$ $e_j,$ $e_k,$ $e_p,$ $e_q,$ $e_r$ along the triangular faces, then
$\Theta(e_i)+\Theta(e_j)+\Theta(e_k)+ \Theta(e_p)+\Theta(e_q)+\Theta(e_r) < 3\pi.$
\item If faces $F_i$ and $F_j$ meet along an edge $e_{ij}$, faces $F_j$ and
$F_k$ meet along an edge $e_{jk}$, and $F_i$ and $F_k$ intersect in exactly one
ideal vertex distinct from the endpoints of $e_{jk}$ and
$e_{ij}$, then $\Theta(e_{ij})+\Theta(e_{jk}) < \pi$.  
\end{enumerate}
Up to isometry, the realization of an abstract polyhedron is unique.  The
ideal vertices of the realization are exactly those degree $3$ vertices for
which there is equality in condition (2) and the degree $4$ vertices.
\end{theorem}

We recall some necessary facts about orbifolds here, and refer the reader to several excellent 
resources \cite{BMP03-1,CoopHodgKer00}. A \emph{hyperbolic Coxeter 
polyhedron} 
is a hyperbolic polyhedron all of whose dihedral angles are integer 
submultiples of $\pi$. To any 
finite volume hyperbolic Coxeter polyhedron $\PP$, there corresponds  a hyperbolic 
$3$--orbifold $\MO_{\PP}$ obtained as the quotient space of $\mathbb{H}^{3}$ 
by the discrete group of isometries generated by all the 
rotations of the form $\rho \sigma$, where $\rho$ and $\sigma$ are reflections 
in two adjacent faces of $\PP$. We call $\MO_{\PP}$ a
\emph{hyperbolic polyhedral $3$--orbifold}. It is topologically the $3$--sphere 
with a marked graph that corresponds to the $1$--skeleton of $\PP$ 
(a general polyhedral $3$--orbifold is topologically 
the $3$--sphere with a marked graph 
that corresponds to the $1$--skeleton of a polyhedron that is not necessarily
hyperbolic). 
A \emph{hyperbolic turnover} is a $2$--dimensional orbifold with underlying space the $2$--sphere
and singular locus consisting of three integer--marked points for which the sum of the reciprocals of the integral markings
is less than $1$.   It is a consequence of a theorem of Dunbar (\cite{Dunbar88-1}, \cite[Theorem 4.8]{BMP03-1}) that the 
$3$--dimensional polyhedral orbifold corresponding to a  hyperbolic Coxeter polyhedron can be decomposed 
(uniquely, up to isotopy) along a system of essential, pairwise non-parallel hyperbolic turnovers into compact, 
irreducible components (with turnover boundary components, if the system is nonempty) 
that contain no essential (embedded) turnovers, 
and such that each component is of one of the following types:
\begin{enumerate}
	\item A $3$--orbifold that contains an essential $2$--suborbifold (that is not a turnover), or
	\item A $3$--orbifold that contains no essential $2$--suborbifolds but that is not the product of a 
	turnover with an interval, or
	\item A $3$--orbifold that is the product of a hyperbolic turnover with an interval.
\end{enumerate}
 This decomposition is equivalent to cutting the planar projection 
  of $\partial P$ along all prismatic $3$--circuits
   (cf. \cite[Corollary 3]{atkinson-thesis}). 
 We call the collection of components of types (1) and (2) the 
 \emph{Dunbar decomposition} of the hyperbolic polyhedron. Components of type (1)
 are called \emph{Haken}, and 
 components of type (2) are called \emph{small}. A hyperbolic Coxeter polyhedron is 
 called \emph{Haken} or \emph{small} if its decomposition consists of a single component 
 of type (1) or (2), respectively. 
  
\begin{figure}
		\centering
		\def\svgwidth{1.5in}
	    	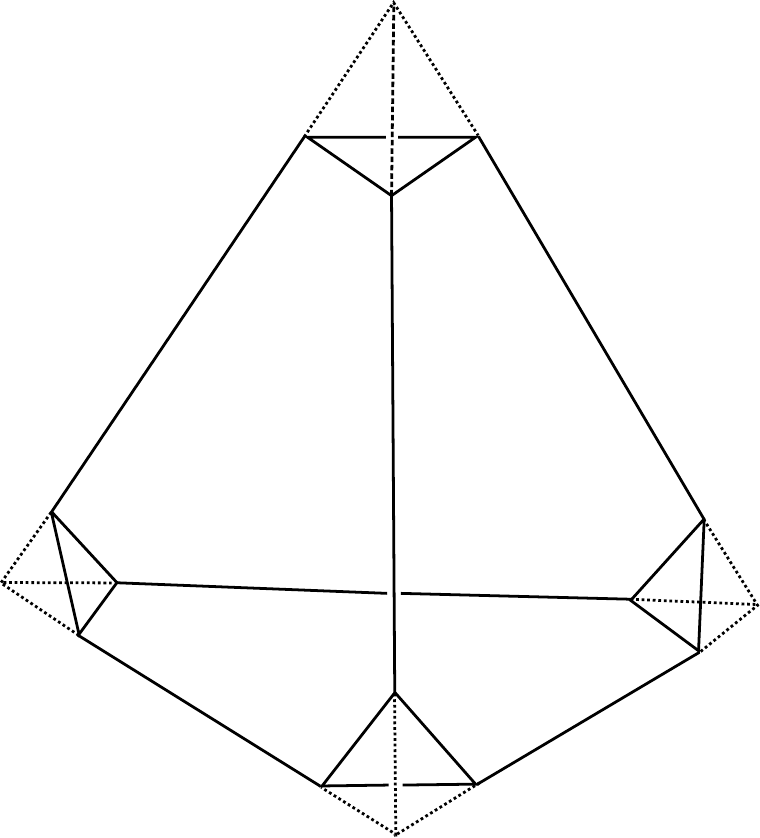
		\caption{The small Coxeter polyhedra in $\mathbb{H}^{3}$}
		\label{F:GenlzdTet}
\end{figure}

A hyperbolic turnover in the $3$--orbifold corresponding to a hyperbolic
polyhedron can  either be made totally
geodesic by an isotopy, or else it doubly covers an embedded, non-orientable
totally geodesic triangular $2$--orbifold that corresponds to a triangular
face, with all right dihedral angles, of the polyhedron (e.g.  \cite[Chapter
IX.C]{Maskit87}, \cite[Theorem 2.1]{Adams05-1}).  Therefore, the Dunbar
decomposition of a hyperbolic polyhedron divides the volume of the polyhedron
additively, and so it is natural, in attempting to organize hyperbolic
polyhedral volumes, to consider the organization up to this decomposition.
The small hyperbolic Coxeter polyhedra have been classified by Rafalski 
 \cite[Theorem 1.1]{rasmallski}:

\begin{theorem}\label{rafalski-small}
A $3$--dimensional hyperbolic Coxeter polyhedron is small if and only
if it is a generalized tetrahedron (see Figure \ref{F:GenlzdTet}).
\end{theorem}

A \emph{generalized hyperbolic tetrahedron} is a hyperbolic tetrahedron 
for which any or all of the vertices are allowed to be ideal or the truncations of hyperideal 
points (see the discussion of \emph{truncated tetrahedra} 
on page \pageref{PR:trunc-simplex}). 
There are nine generalized tetrahedra with all four generalized vertices being finite 
(e.g. \cite[Chapter 7]{Ratcliffe94-1}). 
Every other generalized tetrahedron is obtained from one (or more) of these nine by 
decreasing the dihedral angles of some collection of edges, a process that increases 
volume, by Schl\"{a}fli's formula \cite{milnor}. 
Of the nine generalized tetrahedra with all finite generalized vertices, the one 
of lowest volume is the $3$--$5$--$3$ Coxeter tetrahedron. Its volume  
to six decimal places is $0.039050$. Theorem \ref{T:SmallestLarge} determines 
the lowest volume Haken hyperbolic polyhedron.

\section{On polyhedral volumes}\label{s:volumes}

The main result of this section is Proposition~\ref{graph-volume} in which we
improve on a result of Atkinson giving a lower bound on the
volume of a hyperbolic polyhedron of graph type \cite{atkinson-thesis}.
Theorem \ref{T:SmallestLarge} also requires techniques developed by Atkinson
 \cite{atkinson-equiangular, atkinson-thesis} which we recall here.

We begin with the following theorem that bounds the volume of a hyperbolic
polyhedron without prismatic $4$--circuits in terms of the number of vertices
\cite[Theorem 1.1]{atkinson-thesis}.   This theorem will be used in our proof
of Theorem~\ref{T:SmallestLarge} in Section~\ref{s:proof} to show that any polyhedron with volume less
than that of the Lambert cube must contain a prismatic $4$ circuit.

\begin{theorem}\label{atkinson-lower}
Let $\PP$ be a non--obtuse hyperbolic polyhedron containing
no prismatic $4$--circuits, $N_4$ degree--$4$ vertices, and $N_3$ degree--$3$
vertices. Then
	$$\frac{4N_4 + N_3 -8}{32} \cdot V_8  < 
	\vol(\PP),$$
where $V_8 = 3.663862\dots$ is the volume of the regular ideal hyperbolic
octahedron.
\end{theorem}

The idea behind the proof of Theorem~\ref{atkinson-lower} is to deform $\PP$
to a right-angled polyhedron through an angle-nondecreasing family of
polyhedra. By Schl\"{a}fli's formula, this deformation does not increase volume.  
The assumption that there are no prismatic $4$--circuits ensures
that all polyhedra in this family are hyperbolic. 
The lower bound then comes from applying
Miyamoto's theorem, discussed below, to the resulting right-angled
polyhedron \cite{atkinson-equiangular}.

In order to prove Theorem~\ref{T:SmallestLarge}, we will also need to show
that polyhedra with certain types of prismatic $4$--circuits also have volumes
exceeding that of the Lambert cube.  However,  the techniques used in
Theorem~\ref{atkinson-lower} do not work in the presence of prismatic
$4$--circuits, because attempting such a deformation can cause some or all of
the polyhedron to degenerate to a Seifert fibered polyhedral orbifold.
Although there are volume bounds for such polyhedra \cite[Theorem
13]{atkinson-thesis}, they give a lower bound of $0$ for an infinite family
of polyhedra of graph type. In Proposition~\ref{graph-volume}, we improve
this to give a non--zero lower bound for all hyperbolic polyhedra of graph
type.  A polyhedron $\PP$ is of \emph{graph type} if the polyhedron
$\PP^{\perp}$ obtained by replacing all dihedral angles by $\pi/2$
corresponds to a graph orbifold (i.e., an orbifold with no atoroidal
components in its geometric decomposition).  Atkinson has  shown that 
 all such polyhedra are obtained
by gluing right-angled prisms along quadrilateral faces, and that the components
of the geometric decomposition of the $3$--orbifold for $\PP^{\perp}$
correspond exactly to the prisms that were glued \cite[Sections 4.2 and
6.4]{atkinson-thesis}.  

The remainder of this section focuses on proving
Proposition~\ref{graph-volume}.  The proof is an application of Miyamoto's
theorem, which we recall below.  We first describe the relevant background.

\label{PR:trunc-simplex}
In the projective model of $\HH^3$, thought of as the open ball in $\RR^3
\cup \{\infty\}$, consider any linearly independent set of $4$ points that
lie outside $\HH^3$  such that the line
between any two points meets $\overline{\HH^3}$ in at least one point.  The
intersection of the convex hull of these points with $\HH^3$ is an infinite
volume polyhedron, provided at least one of the points lies outside
$\overline{\HH^3}.$  Form a finite volume polyhedron by truncating the
vertices that lie outside $\overline{\HH^3}$ by their polar hyperplanes.  The
result is a \emph{truncated tetrahedron}.  A truncated tetrahedron is
$\emph{regular}$ if the every edge between a pair of truncating planes has
the same length.

Define $\rho_{3}\left(r\right)$ to be the ratio of the volume of the
regular truncated tetrahedron of edge length $r$ to the sum of the
areas of its faces.  The dihedral
angle $\theta$ along an edge of a regular truncated tetrahedron is
determined by the edge length $r$ via the formula
\begin{equation}\label{E:rVStheta}
\cosh{r} = \frac{\cos{\theta}}{2 \cos{\theta}-1}.
\end{equation}  

Denote the truncated
tetrahedron with dihedral angle $\theta$ by $T_{\theta}$.
The explicit calculation of $\rho_{3}(r)$ was given by Miyamoto \cite{miyamoto}:

\begin{align}\label{E:Rho}
\rho_3(r) &= \frac{\vol(T_{\theta(r)})}{4(\pi - 3 \theta(r))}  \\
          &= \frac{1}{4(\pi - 3 \theta(r))}\left(
		  	V_8 - 3 \int_0^{\theta(r)} \cosh^{-1} \left(\frac
			{\cos (t)}{2\cos (t) - 1}\right) \, dt \right). \notag
\end{align}
where $\theta(r)$ is defined implicitly by Equation \ref{E:rVStheta}.
Miyamoto also proved that $\rho_{3}$ is increasing in  
$r$ \cite[Lemma 2.1]{miyamoto}. Using Equation \ref{E:rVStheta}, $\rho_{3}$ is also 
easily shown to be increasing in $\theta$.

A \emph{return path} in an orbifold with totally geodesic boundary is an
orthogeodesic segment with endpoints on the geodesic boundary.

We can now state Miyamoto's theorem.  Miyamoto proved this theorem for all
dimensions $n\geq 2$, but we require only the $3$--dimensional case.

\begin{proposition}\label{P:Miyamoto1}
	If a complete hyperbolic $3$--orbifold $Q$ of finite volume with totally
	geodesic boundary has a lower bound $l \geq 0$ for the length of its
	return paths, then $$\Vol(Q) \geq \rho_{3}\left({l \over 2}\right)
	\area(\partial Q).$$
\end{proposition}

If $F_{1},...,F_{r}$ is a collection of faces of a hyperbolic polyhedron
$\PP$ such that the sides of each $F_{i}$ all have dihedral angles of $\pi/2$
in $\PP$ and such that no two faces share an edge, then this collection of
faces corresponds to a collection of  embedded, totally geodesic polygonal
$2$--suborbifolds with mirrored sides in the hyperbolic polyhedral
$3$--orbifold $\MO_{\PP}$ obtained from $\PP$. If $Q$ is the hyperbolic
$3$--orbifold with totally geodesic boundary obtained by cutting $\MO_{\PP}$
along this collection of $2$--orbifolds and taking the metric completion,
then Miyamoto's theorem applies. In particular, the inequality that gives the
volume bound may be divided by 2 to give 
$$\Vol(\PP) \geq \rho_{3}\left({l
\over 2}\right) \sum_{i=1}^{r}\Area(F_{i}).$$
 
To get the best volume bound from Miyamoto's theorem, we use the following
special case of a
proposition of Rafalski \cite[Proposition 8.2]{rafalski-thesis} to calculate
the lower bound for the length of a return path in an orbifold in
terms of the Euler characteristic of the boundary.  A return path is
\emph{closed} if its endpoints are equal.  This can happen, for example, if
the orthogeodesic encounters an order--$2$ elliptic axis.

\begin{proposition}\label{P:Miyamoto2}
	Let $Q$ be a complete hyperbolic $3$--orbifold with closed totally
	geodesic boundary.  Then $Q$ has a shortest return path $\gamma$, and
	there is a positive integer $k$ such that the length of $\gamma$ is at
	least the edge length of $T_{\th}$, where 
    $$\th = {\p \over 3(1-k\chi(\partial Q))},$$ 
	where $k>1$ if and only if $\gamma$ is contained in a singular
	axis of (maximal) order $k$ in $Q$ (i.e., $\gamma$ meets a cone point of 
	order $k$ in $\partial Q$). 
\end{proposition}

Combining Proposition~\ref{P:Miyamoto2} with
Proposition~\ref{P:Miyamoto1}, we obtain the following theorem, which gives a lower
bound on the volume of a hyperbolic $3$--orbifold with totally geodesic
boundary in terms of the Euler characteristic of the boundary.  

\begin{theorem}\label{T:Miyamoto3}
Suppose that $Q$ is a hyperbolic $3$--orbifold with totally geodesic boundary
$\partial Q$.  Let 
$$f(t) = \cosh^{-1}\left( \frac{\cos{t}}{2\cos{t} - 1} \right).$$
Then if $k$ is the maximal order of an elliptic element of
$\pi_1(\partial Q)$,  $x=-\chi(\partial Q),$  and 
$$R = 1/2 \cdot  f\left(\frac{\pi}{3(1+kx)}\right),$$
we have
$$\vol(Q) \geq \frac{2 \pi x}{4(\pi - 3 f^{-1}(R))} \left(V_8 -
3\int_0^{f^{-1}(R)} f(t)\, dt \right).$$
\end{theorem}

\begin{proof}
Miyamoto proved that $\rho_3(r)$ is increasing with respect to both $r$ and
$\theta$.  Because the lower bound for $l$ is determined, in Proposition \ref{P:Miyamoto2}, 
as the edge length of a regular truncated tetrahedron $T_{\theta}$ (where 
$\theta$ depends on the Euler characteristic of the orbifold boundary), and
because edge length is minimized when $\theta$ is minimized, we consider the
conditions on a return path that minimize $\theta$. Referring to the
statement of Proposition~\ref{P:Miyamoto2}, we can conclude that $\theta$ is
minimized when $k$ is largest.  The conclusion of the theorem follows by
using Equation~\ref{E:rVStheta} along with Proposition~\ref{P:Miyamoto2} to
expand the volume bound given by Proposition~\ref{P:Miyamoto1}.
\end{proof}

To a polyhedron $\PP$ of graph type that is not a prism, we associate a
graph $G(\PP)$ with vertex set consisting of the set of prisms in the
canonical geometric decomposition of $\PP^{\perp}$ and edges between any two vertices
whose corresponding geometric pieces are glued along a quadrilateral.  A
similar definition could be made in the case of general graph orbifolds.
Note that in the case of polyhedral orbifolds, $G(\PP)$ is a tree, and so has at
least two degree--$1$ vertices.  The number of degree--$1$ vertices of
$G(\PP)$ can be used to bound the volume of $\PP$ below.  The following lemma then
shows that each degree--$1$ vertex of $G(\PP)$ yields a quadrilateral face
that corresponds, after a volume-nonincreasing deformation of $\PP$, to an
embedded totally geodesic suborbifold in the resulting $3$--orbifold.

\begin{lemma}\label{twoquads}
Suppose that $\PP$ is a hyperbolic polyhedron of graph type (that is 
not a prism) such that
$G(\PP)$ has $m$ degree--$1$ vertices.  Then there exists a
volume--nonincreasing deformation of $\PP$ to a hyperbolic polyhedron $\PP'$ 
for which $\MO_{\PP'}$ contains a totally geodesic suborbifold 
of area at least $\frac{m\pi}{6}.$
\end{lemma}

\begin{proof}

We show that there exists an angle
non-decreasing deformation of $\PP$ to a polyhedron $\PP'$ such that all edges of the
quadrilaterals have dihedral angle $\pi/2$, and that $\PP'$ satisfies Andreev's
theorem.

Each degree--$1$ vertex $v$ corresponds to a Seifert fibered component $Q_v$
whose singular locus has the combinatorial type of an $n$--prism, $n \geq 5$.
The fact that $n \geq 5$ follows from the fact that the vertices of $G(\PP)$
correspond to the components of the canonical geometric decomposition.
Degree--$1$ implies that $Q_v$ meets only one other component of the
geometric decomposition, leaving $n-3$ quadrilateral faces free, sharing
edges in a linear fashion.  In each $n$--prism, choose one of these
quadrilaterals.

There exists and angle--increasing
deformation of a polyhedron in which each of the edges along each of the
quadrilaterals has dihedral angle $\pi/2$ \cite[Lemma 12]{atkinson-thesis}.  Then the interior angles of the
quadrilaterals are equal to the dihedral angles along the edges emanating from
the corresponding vertices. Since the area of the quadrilateral is 
the difference of $2\pi$ and its interior angle sum, 
the smallest hyperbolic Coxeter quadrilateral
has area $\pi/6$. 
\end{proof}

The preceding three results can be used to prove the following:

\begin{proposition}\label{graph-volume}
Suppose that $\PP$ is a hyperbolic polyhedron of graph type, and that $\PP$ is not
a prism.  Let $\MC$ be the Lambert cube.  Then

$$\vol(\PP) > \vol(\MC) = 0.324423....$$
\end{proposition}

\begin{proof}
Suppose that $m$ is the
number of degree--$1$ vertices of $G(\PP)$.
We may deform $\PP$ in 
a volume non--increasing manner to a polyhedron $\PP'$, where $\PP'$ has all 
dihedral angles of the form $\pi/2$ and $\pi/3$ \cite[Lemma 12]{atkinson-thesis}.
By Lemma~\ref{twoquads}, there exists a collection of 
$m$ quadrilateral faces of $\PP'$ all of whose edges are labeled 
2 and with total area at least
$\frac{m\pi}{6}$. Each quadrilateral in this collection has between one and four vertices  
with an interior angle equal to $\pi/3$ (corresponding to the dihedral angle between the two 
faces of the polyhedron that meet the quadrilateral at each such vertex). 
The other vertices have interior angles of $\pi/2$. Let $m_{i}$ ($i \in \{1,2,3,4\}$) denote
the number of quadrilateral faces with $i$ interior angles of $\pi/3$. Then $m=\sum_{i} m_{i}$.

By the discussion in the paragraph following 
Proposition \ref{P:Miyamoto1}, we have 
\begin{equation}\label{E:MiyamotoBound}
\Vol(\PP') \geq \frac{\pi}{6} \left(\sum_{i=1}^{4} im_{i} \right) \rho_{3}\left( \frac{l}{2} \right),
\end{equation}
where $l \geq 0$ is a lower bound for the length of a return path in the orbifold with 
boundary obtained by cutting the $3$--orbifold $\MO_{\PP'}$ along the 
collection of totally geodesic suborbifolds corresponding to the $m$ quadrilaterals of $\PP'$. 

We recall that $\rho_{3}$ is increasing with respect to both $r$ and $\theta$ (where $r$ 
is the edge length of the regular truncated tetrahedron $T_{\theta}$ with dihedral angle $\theta$).
In particular, $\rho_{3}(0) = 0.291560...$ is the minimum of this 
function on its domain.    
By Equation \ref{E:rVStheta}, $r$ and $\theta$ each increase 
with respect to the other. Using the minimum of $\rho_{3}$, 
it is easily shown that the lower bound for volume given above is greater than $\Vol(\MC)$ 
if either of $m_{4}$ or $m_{3}$ is nonzero. The bound is also larger than $\Vol(\MC)$ if 
$m_{2} \geq 2$, if $m_{2}=1$ and $m_{1}$ is nonzero, or if $m_{1} \geq 3$.
Since $m \geq 2$, we are left to consider the case when $m_{1} = 2$. 

The minimum of $\rho_{3}$ is insufficient to give the 
appropriate lower bound for $\Vol(\PP')$ in this case.  Using Theorem~\ref{T:Miyamoto3} with $k = 3$ gives a lower bound of $0.406419...$ which is
larger than $\Vol(\MC),$ completing the proof. \end{proof}

\section{Proof of the Main Theorem}\label{s:proof}

The idea behind the proof of Theorem \ref{T:SmallestLarge}
is to use Proposition~\ref{graph-volume} in conjunction with techniques
established by Atkinson \cite{atkinson-equiangular, atkinson-thesis}
and Inoue \cite{inoue} to
restrict the possible combinatorial types of polyhedra with small volumes. 
All polyhedral volumes in the proof are calculated using known formulae and 
the computational software Orb, developed by Heard \cite{HeardOrb}.

\begin{thma}
The smallest volume Haken hyperbolic Coxeter polyhedron is the Lambert
cube $\MC$.
\end{thma}

\begin{proof} 

Suppose that $\PP$ is a Haken hyperbolic Coxeter polyhedron such that
$\vol(\PP) < \vol(\MC) = 0.324423....$  We will show that $\PP$ must have the combinatorial type of a prism.

First note that $\PP$ must contain at least one prismatic $4$--circuit.  If
$\PP$ contained no prismatic $4$--circuits, then $\PP^{\perp}$ admits a
structure as a compact right--angled hyperbolic polyhedron such that
$\vol(\PP^{\perp})\leq \vol(\PP).$  The smallest volume
compact, right--angled polyhedron is the right-angled dodecahedron 
\cite{inoue}. But the volume of this polyhedron is
$4.306207...$, which is considerably larger than $\Vol(\MC)$.  

So we suppose that $\PP$ has at least one prismatic $4$--circuit.  There are
two cases to consider: either the geometric decomposition 
of the orbifold $\MO_{\PP^{\perp}}$ for
$\PP^{\perp}$ contains at least one atoroidal component, or 
$\MO_{\PP^{\perp}}$ is a graph orbifold (i.e., $\PP$ 
is of graph type).

In this former case, suppose that $Q$ is an atoroidal component of 
$\MO_{\PP^{\perp}}$. Then $Q$ corresponds to a 
finite volume right--angled hyperbolic polyhedron $\QQ$ as well as 
to a subset $Q_{\PP}$ of $\PP$ for which we have 
$\vol(\QQ) \leq \vol(Q_{\PP}) \leq \vol(\PP)$ \cite[Propositions 2 and 3]{atkinson-thesis}.
The proof of this fact uses a result of Agol, Storm and Thurston 
\cite[Corollary 2.2]{agol-storm-thurston}. 
The polyhedron $\QQ$ has at least one $4$--valent ideal vertex
coming from the prismatic $4$--circuit.  The following lemma shows that the
volume of $\QQ$  must exceed that of the Lambert cube.

\begin{lemma}\label{no2vtxrtanglepoly}
Suppose $\QQ$ is a right--angled hyperbolic polyhedron with at least one ideal
vertex.  Then
$\vol(\QQ) > \vol(\MC)$.
\end{lemma}
\begin{proof}
Let $N_3$ and $N_4$ be the number of $3$--valent and $4$--valent vertices of
$\QQ$ respectively.  If $N_4 \geq 3$, then Theorem~\ref{atkinson-lower}
immediately gives the conclusion in this case because 
$$\vol(\QQ) \geq \frac{4N_4 -8}{32}\cdot V_8 > 0.457 > \vol(\MC).$$

We now discuss the cases where $N_4$ is either $1$ or $2$. If $E$ is the
number of edges of $\QQ$, note that $2E = 3N_3 + 4N_4.$  

Suppose first that $N_4 = 1$.  Then $2E = 3N_3 + 4$, so $N_3$ is even.  The minimal
number of vertices of a polyhedron is $4$, so $N_3 \geq 3$.  The only example
with $N_4 =1$ and $N_3 = 4$ is an ideal pyramid with square base which, by
Andreev's theorem,  does not admit a right angled hyperbolic realization.
This discussion along with the following lemma shows that a right--angled
hyperbolic polyhedron with a single ideal vertex must have more than $6$
finite vertices.

\begin{lemma}\label{no1vtxrtanglepoly}
There are no right-angled hyperbolic polyhedron with $N_4 = 1$ and $N_3 = 6$.
\end{lemma}

\begin{proof}
Suppose that $\PP$ is such a polyhedron.  Label the ideal vertex $v_0$, the
four vertices adjacent to $v_0$ by $v_1$, $v_2$, $v_3$, and $v_4$, and the
remaining two vertices by $v_5$ and $v_6$.  By Euler characteristic
considerations, $\PP$ must contain $11$ edges.  Denote the set $\{v_1, v_2,
v_3, v_4\}$ by $\VV$.  We first observe that $\PP$ must contain at least one
edge with both endpoints in $\VV$.  To see this, note that each vertex in
$\VV$ is trivalent, so each vertex in $\VV$ meets two more edges in addition
to the edge shared with $v_0$.  There are only seven edges in addition to the
edges containing $v_0$, so by the pigeonhole principle, at least two of the
vertices in $\VV$ must share an edge.  

In fact, there must be two edges with both endpoints in $\VV$.  Suppose that
there was only one such edge.  We may assume, without loss of generality,
that the edge with both endpoints in $\VV$ meets $v_1$ and $v_2$.  Vertices
$v_3$ and $v_4$ must have valence $3$, so both
must meet two additional edges.  Bigons are not permitted, so each of $v_3$ and $v_4$ must
meet edges containing $v_5$ and $v_6$ giving a configuration as in
Figure~\ref{two-edges}, up to switching $v_5$ and $v_6$.  Suppose for
concreteness that the embedding of the $1$--skeleton is as shown.  The
edge path of length four passing through $v_0$, $v_3$, $v_6$, and $v_4$
separates $v_5$ from $v_1$ and $v_2$.  This is a contradiction to planarity
because $v_5$ must be connected to one of $v_1$ and $v_2$ by an edge.

\begin{figure}
\labellist
\small \hair 2pt
\pinlabel {$v_1$} [br] at 1 65
\pinlabel {$v_2$} [tr] at 1 1
\pinlabel {$v_3$} [br] at 64 65
\pinlabel {$v_4$} [tr] at 64 1
\pinlabel {$v_5$} [l] at 84 32
\pinlabel {$v_6$} [l] at 115 32
\pinlabel {$v_0$} [l] at 35 32
\endlabellist
\scalebox{1}{\includegraphics{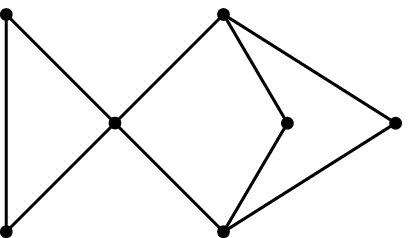}}
\caption{}
\label{two-edges}
\end{figure}

Finally, we observe that there cannot be $3$ or $4$ edges with both endpoints
in $\VV$.  Having $3$ such edges forces a bigon between $v_5$ and $v_6$ and
$4$ such edges causes the graph to be disconnected.

The three possibilities for a such a polyhedron with two edges having both
endpoints in $\VV$  are shown in Figure~\ref{poly-cases}.  In the first two cases, the graphs are not $3$--connected.
The third case contradicts the fourth condition of Andreev's theorem.
\end{proof}

\begin{figure}
\scalebox{.7}{\includegraphics{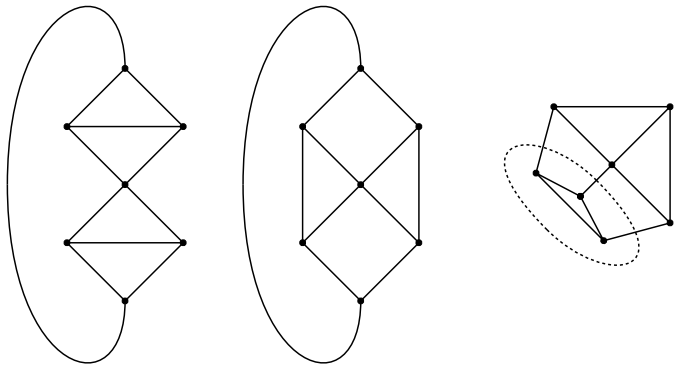}}
\caption{The first two graphs are not $3$--connected.  The dashed curve in
the third graph illustrates the contradiction to Andreev's theorem.}
\label{poly-cases}
\end{figure}

It follows from Lemma~\ref{no1vtxrtanglepoly} and discussion 
preceding it that a right--angled hyperbolic
polyhedron with a single ideal vertex must have at least $8$ finite vertices.
Theorem~\ref{atkinson-lower} then says 
$$\vol(\QQ) > \frac{4\cdot 1 + 8 - 8}{32} \cdot V_8 > 0.458 > \vol(\MC).$$ 

If $N_4 = 2$, then $2E = 3N_3 +8$, so $N_3$ is odd.  Then since
$N_3 \geq 3$  Theorem~\ref{atkinson-lower}  gives
$$ \vol(\QQ) > \frac{ 4 \cdot 2 + 3 - 8}{32} \cdot V_8 > 0.343 >
\vol(\MC).$$

This completes the proof of the Lemma~\ref{no2vtxrtanglepoly}.
\end{proof}

To complete the proof of Theorem~\ref{T:SmallestLarge}, 
we are left to consider the case when $\PP$ is
a polyhedron of graph type.  If $\PP$ is not a prism, then
Proposition~\ref{graph-volume} implies that the volume of $\PP$
exceeds that of $\MC$.  Hence if $\vol(\PP)\leq\vol(\MC),$ $\PP$ must have the
combinatorial type of a prism.

Any hyperbolic Coxeter $3$--prism (i.e., a triangular prism) is either 
a generalized tetrahedron with one truncated vertex 
or has a Dunbar decomposition into two generalized tetrahedra each 
with one truncated vertex \cite[Lemma 3.1]{rasmallski}. 
In either case, such a prism is not Haken.
Atkinson has determined the smallest volume Coxeter $n$--prisms for $n\geq 5$, 
and shown that the smallest volumes increase monotonically in $n$ 
\cite[Theorem 11 and Corollary 8]{atkinson-thesis}. 
The lowest $5$--prism volume is $0.763304...$, greater than that of $\MC$.  
The $4$--prism case remains, that is, we must show
that $\MC$ has the smallest volume among Coxeter polyhedra with the
combinatorial type of the cube.

For any Coxeter polyhedron with the combinatorial type of the cube, there
exists an volume nonincreasing deformation to one with all dihedral angles
$\pi/2$ and $\pi/3$ \cite[Proposition 4]{atkinson-thesis}.  The only
restriction that Andreev's theorem places on such a polyhedron is that there
is at least one dihedral angle of $\pi/3$ along each of the three prismatic
$4$--circuits.  If there were more than one dihedral angle of $\pi/3$ along
one of the prismatic $4$--circuits, there would exist a further volume
decreasing deformation to a polyhedron with only one dihedral angle of
$\pi/3$.  Hence any Coxeter polyhedron with the combinatorial type of the
cube has volume at least that of a cube with a single dihedral angle of
$\pi/3$ occurring on each of the prismatic $4$--circuits and all other
dihedral angles $\pi/2$.  Up to isometry, there are $4$ such polyhedra,
$\MC_1 = \MC$, $\MC_2$, $\MC_3$, and $\MC_4,$ as
shown in Figure~\ref{cubes}. 
The respective volumes of   $\MC_1$, $\MC_2$,
$\MC_3$, and $\MC_4,$ are $0.324423\dots$, $0.392365\dots$, $0.464467\dots$,
and $0.634337\dots$.

\begin{figure}
\labellist
\small\hair 2pt
\pinlabel {$\frac{\pi}{3}$} [r] at 0 33
\pinlabel {$\frac{\pi}{3}$} [b] at 34 46
\pinlabel {$\frac{\pi}{3}$} [tr] at 55 14
\pinlabel {$\frac{\pi}{3}$} [r] at 92 33
\pinlabel {$\frac{\pi}{3}$} [b] at 126 21
\pinlabel {$\frac{\pi}{3}$} [tr] at 148 14
\pinlabel {$\frac{\pi}{3}$} [r] at 204 33
\pinlabel {$\frac{\pi}{3}$} [b] at 218 21
\pinlabel {$\frac{\pi}{3}$} [tr] at 240 14
\pinlabel {$\frac{\pi}{3}$} [r] at 297 33
\pinlabel {$\frac{\pi}{3}$} [b] at 309 21
\pinlabel {$\frac{\pi}{3}$} [tl] at 288 14
\pinlabel {$\MC_1$} [t] at 35 -1
\pinlabel {$\MC_2$} [t] at 128 -1
\pinlabel {$\MC_3$} [t] at 219 -1
\pinlabel {$\MC_4$} [t] at 310 -1
\endlabellist
\scalebox{1}{\includegraphics{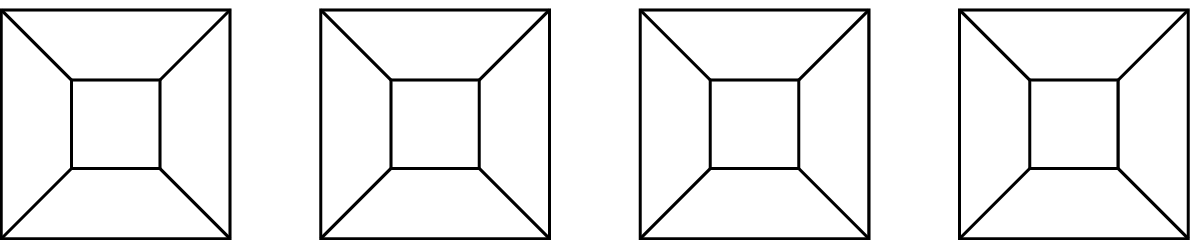}}
\caption{Lowest volume candidates for combinatorial cubes}
\label{cubes}
\end{figure}

This completes the proof.
\end{proof}

\bibliographystyle{hamsplain}
\bibliography{AtkinsonRafalski_refs}

\end{document}